\documentclass[a4paper,onecolumn,oneside,11pt]{article}
\usepackage{color}
\usepackage{amsmath, amsfonts,amssymb,theorem,
euscript,array,enumerate,amsfonts,mathrsfs}
\usepackage[dvips]{graphicx}
\usepackage[dvips]{graphics}
\usepackage[dvips]{epsfig}

\newcommand{\1}{\mathbf{1}}
\newcommand{\ind}{\mathbf{1}}

\newcommand{\cgblue}{\color{blue}}

\definecolor{ying}{rgb}{0.3,0.1,0.8}
\theoremheaderfont{\bfseries} \theoremstyle{plain}
\theorembodyfont{\upshape}
\newcommand{\esp}{\mathbb{E}} 
\newcommand{\E}{\mathbb{E}}
\newcommand{\indic}{1\!\!\!1} 
\newcommand{\proba}{\mathbb{P}}

\newcommand{\ff}{\mathbb{F}}

\newcommand{\Q}{\mathbb{Q}}
\newcommand{\R}{\mathbb{R}}
\newcommand{\D}{\mathcal{D}}
\newcommand{\F}{\mathcal{F}}
\newcommand{\G}{\mathcal{G}}

\newcommand{\cA}{\mathcal{A}}

22 truecm \hbadness 10000 \oddsidemargin=-0,5cm \tolerance 10000

\theoremstyle{plain}
\theorembodyfont{\upshape\it}
\newtheorem{Thm}{\bf Theorem}[section]
\newtheorem{Pro}[Thm]{\bf Proposition}
\newtheorem{Lem}[Thm]{\bf Lemma}
\newtheorem{Cor}[Thm]{\bf Corollary}
\theorembodyfont{\upshape}

\newtheorem{Def}[Thm]{Definition}

\newtheorem{Rem}[Thm]{Remark}
\newtheorem{Rems}[Thm]{Remarks}
\newtheorem{Hyp}[Thm]{Hypothesis}

\newcommand {\proof} {\noindent {\sc Proof:  }}
\newcommand {\bproof} {\noindent {\sc Proof:  }}
\newcommand {\finproof} {\hfill $\Box$ \vskip 5 pt }

\begin{document}
\thispagestyle{empty}
\title{What happens after a default: the conditional density approach
\footnote{This work has benefited from financial support by  La Fondation du Risque et F\'ed\'eration Bancaire Fran\c caise. } }
\author{
Nicole El Karoui\thanks{LPMA, Universit\'e Paris VI and CMAP,
Ecole Polytechnique,  email: elkaroui@cmapx.polytechnique.fr} 
\quad Monique Jeanblanc\thanks{D\'epartement de Math\'ematiques, Universit\'e
d'Evry and Institut Europlace de Finance, email: monique.jeanblanc@univ-evry.fr}
\quad Ying Jiao\thanks
{Laboratoire de Probabilit\' es et Mod\`eles Al\' eatoires,  Universit\'e Paris VII, Case 7012, 2 Place Jussieu, 75251 Paris Cedex 05, email: jiao@math.jussieu.fr} }
\date{\today}
\maketitle

\bigskip
 \begin{abstract}
 We present a general model for default time, making precise the
 role of the intensity process, and showing that this process
 allows for a  knowledge of the conditional distribution of the default
 only ``before the default". This lack of  information is crucial
 while working in a multi-default setting. In a single default case,
the knowledge of the intensity process does not allow to compute
the price of defaultable claims, except in the case where
immersion property  is satisfied. We propose in this paper the density approach for default time. The density process will give a
full characterization  of the links between the default time and
the reference filtration, in particular ``after the default time".
We also investigate the description of martingales in the full
filtration in terms of martingales in the reference filtration,
and the impact of Girsanov transformation on the density and intensity processes, and also on the immersion property.

 \end{abstract}
\section{Introduction }
Modelling  default time for a single credit  event has been
largely studied in the literature,   the  main approaches being
the structural, the reduced-form  and the intensity ones. In this
context, most works are concentrated (for pricing purpose) on the
computation of conditional expectation of payoffs, given that the
default has not occurred, in the case where immersion property is
satisfied. In this paper, we are also interested in what happens
after a default occurs: we find it important to investigate the
impact of a default event on the rest of the market and what goes
on afterwards.

Furthermore, in a multi-default setting, it will be
important  to compute the prices of a portfolio derivative on the
disjoint sets before the first default, after the first and
before the second and so on. Our work will allow us to use a
recurrence procedure to provide these computations, which will be
presented in a companion paper \cite{EJJ2009}.

We   start with the knowledge of  the conditional distribution of
the default time $\tau$, with respect to a reference filtration
$\mathbb F=(\F_t)_{t\geq 0}$ and we assume that this conditional distribution
admits a density (see the first section below for a precise
definition). We firstly reformulate the classical computation
result of   conditional expectations with respect to the
observation $\sigma$-algebra $\G_t:= \F_t \vee \sigma (\tau \wedge
t)$ before the default time $\tau$, i.e., on the set $\{t<\tau \}$. The main
purpose is then to deduce what happens after the default occurs,
i.e., on the set $\{\tau\leq t\}$. We shall emphasize that the
density approach is
 suitable in this after-default study and explain why the
intensity approach is inadequate for this case. We present
computation results of $\mathbb G=(\G_t)_{t\geq 0}$ conditional expectations on the
set $\{\tau\leq t\}$ by using the conditional density of $\tau$
and point out that the whole term structure of the density is
needed. By establishing an explicit link between (part of) density
and intensity, which correspond respectively to the additive and
multiplicative decomposition related to the survival process
(Az\'ema supermartingale of $\tau$), we make clear  that the
intensity can be deduced from the density, but that the reverse
does not hold, except when certain strong assumption, as the
H-hypothesis, holds.

Note that, even if   the ``density'' point of view is inspired by
the enlargement of filtration theory,  we shall not use classical
results on the progressive enlargement of filtration. In fact, we
take the opposite point of view: we are interested in $\mathbb
G$-martingales and their characterization  in terms of $\mathbb
F$-(local) martingales. Moreover, these characterization results
allow us to give a   proof of a  decomposition of $\mathbb
F$-(local) martingales in terms of $\mathbb G$-semimartingales.

We   study how the parameters of the default (i.e., the survival
process, the intensity, the density) are modified by a change of
probability in a general setting (we do not assume that we are
working in a Brownian filtration, except for some examples), and
we characterize changes of probability that do not affect the
intensity process. We pay attention to  the specific case where
the dynamics of underlying default-free processes are changed only
after the default.

 The paper is organized as follows. We first introduce in
Section 2 the different types of information we are dealing with
and the key hypothesis of density. In Section 3, we establish
results on computation of conditional expectations, on the
``before default'' and ``after default'' sets. The H-hypothesis is
then discussed. The dynamic properties of the density process are
presented in Section 4 where we make precise the links between
this density process and the intensity process. In the last section, we present
the characterization of $\mathbb G$-martingales in terms of
$\ff$-local martingales.  We give a Girsanov type property and
discuss  the stability of immersion property and invariance of
intensity.
\section{The Different  Sources of Information}

In this section, we specify the link between the two filtrations
$\mathbb F$ and $\mathbb G$, and make some hypotheses on the
default time.
Our aim is to measure the consequence of a default event in terms of pricing  various contingent claims.\\
We  start as usual with a filtered probability space $(\Omega,
\cA, \mathbb F, \mathbb P)$. Before the default  time $\tau$,
i.e., on the set $\{t<\tau\}$, the $\sigma$-algebra $\F_t$
represents the information accessible to the investors at time
$t$. When the default occurs, the investors will add this new
information (i.e., the knowledge of $\tau$) to  the
$\sigma$-algebra $\F_t$.
\\
More precisely, a strictly positive and finite random variable
$\tau$  (the default time) is given on the probability space
$(\Omega, \cA,   \mathbb P)$. This space is supposed to be endowed
with a  filtration $\mathbb F=(\mathcal F_t)_{t\geq 0}$ which
satisfies the usual conditions, that is, the filtration $\mathbb F
$ is right-continuous and $\mathcal F_0$ contains all
$\proba$-null sets of $\cA$.
\\
One of our goals is to show how the information contained in the
reference filtration $\ff$  can be used to obtain information on
 the distribution of $\tau$. We  assume that, for any $t$, the conditional distribution of $\tau$ with respect to
$\F_t$ is smooth, i.e., that the $\F_t$-conditional distribution
of $\tau$ admits a density
 with respect to some  positive $\sigma$-finite
measure $\eta$ on $\mathbb R^+$. As an immediate consequence, the
unconditional distribution of $\tau$ is absolutely continuous
w.r.t. $\eta$. Another consequence is that $\tau$ can not be an $\ff$-stopping time.\\ 
 In other terms, we introduce the following
hypothesis\footnote{This hypothesis has been discussed by Jacod
\cite{Ja1987} in the initial enlargement of filtration framework.
The same assumption also appears in the dynamic Bayesian framework
\cite{GVV2006}. We do not assume that $\eta$ is finite, allowing
for the specific case of the Lebesgue measure.}, that we call
density hypothesis. This hypothesis will be in force in the paper.
 \begin{Hyp}\label{Hyp:hyp(J)}
 ({\bf Density hypothesis.}) We assume that $\eta$ is  a  non-negative non-atomic measure  on
$\mathbb R^+$ such that, for any time $t\geq 0$, there exists an
$\F_t\otimes \mathcal B(\mathbb R^+)$-measurable function
$(\omega, \theta) \rightarrow \alpha_t(\omega, \theta)$   which
satisfies
\begin{equation}\label{equ: density def}
\proba(\tau\in d\theta|\F_t)=:\alpha_t(\theta)\eta(d\theta) ,
\,\quad \proba-a.s.\end{equation}
\end{Hyp}
The family $\alpha_t(.)$ is called the {\it conditional density}
of $\tau$ with respect to $\eta$ given $\F_t$  (in short the
{\it density} of $\tau$ if no ambiguity).  Then, the distribution
of $\tau$ is given by ${\mathbb {P}}(\tau >\theta)=\int_\theta
^\infty \alpha _0(u) \eta(du)$. Note that, from the definition and
the hypothesis that $\tau$ is finite, for any $t$, $\int_0^\infty
\alpha_t(\theta)\eta(d\theta)=1\, (a.s.)$. By definition of the
conditional expectation, for any (bounded) Borel function $f$,
$\esp[f(\tau) |\F_t]
 = \int_0^\infty f(u) \alpha _t(u) \eta (du) (a.s.)$.
The conditional distribution of $\tau$ is also characterized by
the survival probability function
\begin{equation}S_t(\theta):= \mathbb P (\tau >\theta \vert \F_t)= \int
_\theta^\infty \alpha_t(u)\eta(du)\end{equation}
 The family of random variables
 $$S_t := S_t(t)=\proba(\tau>t|\F_t)= \int_t^\infty \alpha_t(u)\eta(du)$$ plays a key role in
what follows. Observe that one has
\begin{equation}\label{setA}\{\tau>t\} \subset \{S_t> 0\} =: A_t\end{equation}
(where the inclusion is up to a negligible set)     since
$\proba(A_t^c\cap\{\tau>t\})= 0$. Note also that $S_t(\theta)= \E(
S_\theta \vert \F_t)$ for any $\theta \geq t$.
\\
More generally,  if an $\F_t\otimes \mathcal B(\mathbb
R^+)$-measurable function $(\omega, \theta) \rightarrow
Y_t(\omega, \theta)$  is given, the   $\F_t$-conditional
expectation of the r.v. $Y_t(\tau):= Y_t(\omega, \tau (\omega))$,
assumed to be integrable, is given by
\begin{equation}\label{YTtau}
\E[ Y_t(\tau)\vert \F_t] = \int_0^\infty Y_t(u)\alpha_t(u)\eta(du)
\,.
\end{equation}

\noindent \textbf{Notation:}  In what follows, we shall simply say
that $Y_t(\theta)$ is an $ \F_t\otimes\mathcal B(\mathbb R^+)
$-random variable and even that $Y_t(\tau)$ is an $
\F_t\otimes\sigma(\tau) $-random variable as a short cut for
$Y_t(\theta)$ is an $\F_t\otimes \mathcal B(\mathbb
R^+)$-measurable function.


  \begin{Cor} \label{Cor:avoid F stopping time}
  The default time $\tau$ avoids $\mathbb F$-stopping
  times, i.e., $\proba (\tau = \xi)=0$ for every $\mathbb
  F$-stopping time $\xi$.
  \end{Cor}
  \proof
   Let $\xi$ be  an $\mathbb F$-stopping time bounded  by a
  constant $T$. Then, the random variable $H_\xi(t)=\indic_{\{\xi= t\}}$ is $ \F_T\otimes\mathcal B(\mathbb R^+)
$-measurable, and, $\eta$ being non-atomic \[\E[ H_\xi(\tau)\vert
\F_t] =\E[\,\E[ H_\xi(\tau) \vert \F_T]\,\vert \F_t] =
\E[\int_0^\infty H_\xi(u)\alpha_T(u)\eta(du) \vert \F_t]=0\,.\]
Hence,
  $\E [H_\xi(\tau)]= \proba (\xi=\tau)=0$.
  \finproof

\begin{Rem} \label{forward rate}By using density, we adopt an additive point of view to
represent the conditional probability of $\tau$:       the
conditional survival function
$S_t(\theta)=\proba(\tau>\theta\,|\,\F_t)$ is written on the form
$S_t(\theta)=\int_{\theta}^\infty \alpha_t(u)\eta(du)$. In the
default framework, the ``intensity'' point of view is often
preferred, and  one uses the multiplicative representation
$S_t(\theta)=\exp(-\int_0^{\theta}\lambda_t(u)\eta(du))$. The
family of $\F_t$-measurable random variables
$\lambda_t(u)=-\partial_u\ln S_t(u)$ is called the ``forward
hazard rate". We shall discuss and compare these two points of
view further on.

\end{Rem}

\section{Computation of conditional expectations in a default setting}

The specific information related to the default time is the
knowledge of this time when it occurs. It is defined in
mathematical terms  as follows: let $\mathbb D=(\D_t)_{t\geq 0 }$
be the smallest right-continuous filtration such that $\tau$ is a
$\mathbb D$-stopping time; in other words, $\D_t$ is given by
$\D_t= \D_{t+}^0$ where $\D_{t }^0=\sigma(\tau\wedge t)$. This
filtration $\mathbb D$ represents the default information, that
will be ``added" to the reference filtration; the filtration
$\mathbb G:= \mathbb F\vee \mathbb D$ is  the smallest filtration
containing $\ff$ and making $\tau$ a stopping time. Moreover, any
$\G_t$-measurable r.v. $H^{\mathbb G}_t$ may be represented as
$H^{\mathbb G}_t=H_t^{\mathbb
F}\indic_{\{\tau>t\}}+H_t(\tau)\indic_{\{\tau\leq t\}}$ where
$H_t^{\mathbb F}$ is an $\F_t$-measurable  random variable and $
H_t(\tau) $
 is  $\F_t\otimes \sigma (\tau)$-measurable. In particular,
\begin{equation}\label{della}
H_t^{\mathbb G} \,\indic_{\{\tau>t\}}=H_t^{\mathbb
F}\indic_{\{\tau>t\}} \quad a.s.\,,
\end{equation}
where  the random variable $H_t^{\mathbb F}$     is the
$\F_t$-conditional expectation of $H_t^{\mathbb G}$ given the
event $\{\tau>t\}$, i.e.,
\begin{equation}\label{Equ:cal G cond exp}
H_t^{\mathbb F} =\frac{\esp[H_t^{\mathbb G}
\indic_{\{\tau>t\}}|\F_t]}{\proba(\tau>t|\F_t)}=\frac{\esp[H_t^{\mathbb
G} \indic_{\{\tau>t\}}|\F_t]}{S_t}\quad a.s.\;\mbox{on } A_t;\quad
 H_t^{\mathbb
F}= 0\quad \mbox{on~the~ complementary~set}\,.
\end{equation}
  \subsection{Conditional expectations}
  The definition of $\mathbb G$  allows us to compute
conditional expectations with respect to ${\cal {G}}_t$ in terms
of conditional expectations with respect to ${\cal {F}}_t$. This
will be done   in two steps, depending whether or not the default
has occurred: as we explained above, before the default, the only
information contained in ${\cal {G}}_t$ is  ${\cal {F}}_t$, after
the default, the information contained in ${\cal {G}}_t$ is,
roughly speaking, ${\cal {F}}_t \vee \sigma (\tau)$.

 The
$\G_t$-conditional expectation of an integrable $\sigma
(\tau)$-measurable r.v. (of the form $f(\tau)$) is given by
$$\alpha^{\mathbb G}_t (f):=\esp[f(\tau)|\G_t]=\alpha^{\text{bd}}_t(f)\,\indic_{\{\tau>t\}}
      +f(\tau)\,\indic_{\{\tau\leq t\}}$$ where   $\alpha^{\text{bd}}_t$  is the value of the
$\G_t$-conditional distribution \textbf{b}efore the
\textbf{d}efault, given by
\[\alpha_t^{\text{bd}}(f):=\frac{\esp[f(\tau)\indic_{\{\tau>t\}}|\F_t]}
{\proba(\tau>t|\F_t)} \quad a.s.\,\text{on } \,A_t;
  \quad \alpha_t^{\text{bd}}(f):= 0\quad \mbox{on~the~ complementary~set}. \]

Recall the notation $S_t=\proba(\tau>t|\F_t)$. On the set $A_t$,
the ``before the default" conditional distribution
$\alpha_t^{\text{bd}}$  admits a density $\alpha_t^{\text{bd}}(u)$
with respect to $\eta$, given by
\[\alpha_t^{\text{bd}}(u) =\frac{1 }{S_t }
\ind_{[t,\infty)}(u)\alpha_t(u)\eta(du)  \quad a.s.\, .  \] The
same calculation as in \eqref{YTtau} can be performed in this new
framework and extended to the computation of $\G_t$-conditional
expectations for a
 bounded  $ \F_T \otimes\sigma(\tau) $-r.v..

 \begin{Thm}\label{Thm:cond expect general case with density}
Let $Y_T(\tau)$ be a bounded  $ \F_T \otimes\sigma(\tau) $-random
variable. Then, for $t\leq T$
$$ E[Y_T(\tau) |\G_t] =Y^{\rm{bd}}_t \indic_{\{t<\tau\}}+ Y^{\rm{ad}}_t(T,\tau)\indic_{\{\tau\leq t\}} \quad
d\proba-a.s.
$$
 where
 \begin{eqnarray} Y^{\rm{bd}}_t&=&\frac{\esp \big[\int_ t^\infty Y_T(u)\alpha_T(u)\eta(du)\vert \F_t] }{S_t}
 \nonumber \quad
d\proba-a.s.\, \mbox {on }\,A_t,
\\Y^{\rm{ad}}_t(T,\theta)&=&\frac{\esp\big[Y_T(\theta)\alpha_T(\theta)\big|\F_t\big]}{\alpha_t(\theta)}\,
\indic_{\{\alpha_t(\theta)>0\}} \label{after}\quad d\proba-a.s..
\end{eqnarray}
\end{Thm}
\proof The computation on the set $\{t<\tau\}$ (the pre-default
case) is obtained following (\ref{della}),
(\ref{Equ:cal G cond exp}) and using  (\ref{YTtau}). For the after-default case, we note
that, by definition of $\mathbb G$, any $\G_t$-measurable r.v. can
be written on the set $\{\tau\leq t\}$ as
$H_t(\tau)\indic_{\{\tau\leq t\}}$. Assuming that $H_t(\tau)$ is
positive or bounded, and using the density $\alpha_t(\theta)$, we
obtain
\[\begin{split}
\esp[H_t(\tau)\indic_{\{\tau\leq t\}}Y_T(\tau)] &=\int_0^\infty
d\theta\>\esp[H_t(\theta)\indic_{\{\theta\leq
t\}}Y_T(\theta)\alpha_T(\theta)] =\int_0^\infty
d\theta\>\esp\big[H_t(\theta)\indic_{\{\theta\leq
t\}}\esp[Y_T(\theta)\alpha_T(\theta)|\F_t]\big]\\ &=\int_0^\infty
d\theta\>\esp\Big[H_t(\theta)\indic_{\{\theta\leq
t\}}Y^{ad}_t(T,\theta)\alpha_t(\theta)\Big]=
\esp\Big[H_t(\tau)\indic_{\{\tau\leq t\}}Y^{ad}_t(T,\tau)\Big],
\end{split}
\]
which implies immediately \eqref{after}. \finproof

\subsection{Immersion property or H-hypothesis} \label{subsubsec:(H)-Hypothesis}
In the form of the density $\alpha_t(\theta)=\proba(\tau \in
d\theta \vert \F_t)/\eta (d\theta)$, the parameter $\theta$ plays
the role of the default time. It is hence natural to consider the
particular case where
\begin{equation}\label{Equ:alpha under H-hypothesis}
\alpha_t(\theta)=\alpha_\theta(\theta), \quad \forall \,
\theta\leq t\,,\end{equation} i.e., the case when    the
information contained in the reference filtration after the
default time does not give new information on the conditional
distribution of the default. In that case

$$S_t=
\proba(\tau>t |\F_t) =1-\int_ 0^t\alpha _t(u)\eta (du)=  1-\int_
0^t\alpha _u(u)\eta (du)$$ and, in particular $S$ is decreasing
\footnote{\label{scontinu}and continuous, this last property will
be useful later.}. Furthermore,
$$S_t= 1-\int_ 0^t \alpha _T(u)\eta (du)=\proba(\tau>t
|\F_T)=S_T(t)\; a.s.$$ for any $T\geq t$ and it follows that $
\proba(\tau>t |\F_t) =\proba(\tau>t |\F_\infty)$. This last
equality is known to be equivalent to the immersion property
(\cite{BY1978}), also known as the H-hypothesis, stated as: for any
fixed $t$ and any bounded $\G_t$-measurable r.v. $Y^{\mathbb
G}_t$,
\begin{equation} \label{eq:(H)hyp}
\esp[Y^{\mathbb G}_t|\F_\infty]=\esp[Y^{\mathbb G}_t|\F_t] \quad
a.s..
\end{equation}
Conversely, if immersion property holds, then (\ref{Equ:alpha
under H-hypothesis}) holds. In that case,
 the conditional survival
functions $S_t(\theta)$ are constant in time on $[\theta,
\infty)$, i.e., $S_t(\theta)=S_\theta(\theta)$ for $t>\theta $.
Then the previous result (\ref{after}) takes a simpler form:
$Y_t^{\rm ad}(T,\theta)=\esp[Y_T(\theta)|\F_t],\>a.s.$ for $\theta\leq
t\leq T$, on the set $\{\alpha_{\theta}(\theta)>0\}$. \\
Under immersion property, the knowledge of $S $ implies that of
the conditional distribution of $\tau$ for all positive $t$ and
$\theta$: indeed, one has   $S_t(\theta)=\esp[S_\theta|\F_t]$
(note that, for $\theta <t$, this equality reduces   to
$S_t(\theta)=S_\theta(\theta)=S_\theta$).

\begin{Rems} \label{pseudo} ~\\ 1)  The decreasing property of $S$ (equivalent to the fact
that $\tau$ is a pseudo-stopping time (see \cite{NY2005})) does
not imply the H-hypothesis, but only that  $\ff$-bounded
martingales stopped at $\tau$ are $\mathbb G$-martingales (see
also \cite{EJY2000}). We shall revisit this property
  in Remarks \ref{remaimmer} and \ref{immer} and
Corollary \ref{Cor:martingale continue at tau}.
\\
2) The most important example where immersion holds is the widely
studied Cox-process model introduced by Lando
\cite{La1998}.\end{Rems}


\section{Dynamic point of view and density process}


Our aim is here to give a dynamic study of the previous results.
 We shall call  $(S_t,t\geq 0)$ the survival process, which is
an $\ff$-supermartingale.  We have obtained equalities for fixed
$t$, we would like to study the conditional expectations as
processes. One of the goals is to recover the value of the
intensity of the random time, and the decompositions  of $S$.
Another one is to study the link between $\mathbb G$ and $\ff$
martingales: this is of main interest for pricing.

\noindent In this section, we present the dynamic version of the
previous results in terms of $\mathbb F$ or $\mathbb G$
martingales or supermartingales. To be more precise, we need some
``universal" regularity on the paths of the density process. We
shall treat some technical problems in Subsection \ref{Subsec:path
regularity} which can be skipped for the first reading.

\subsection{Regular Version of Martingales}\label{Subsec:path regularity}
 One of the major difficulties is to
prove the existence of a universal c\`adl\`ag martingale version
of this family of densities. Fortunately, results of Jacod
\cite{Ja1987} or Stricker and Yor \cite{SY1978} help us to solve
this technical problem.

Jacod (\cite{Ja1987}, Lemme 1.8) establishes the existence of a
universal c\`adl\`ag version of the density process in the
following sense: there exists a non negative function
$\alpha_t(\omega,\theta)$ c\`adl\`ag in $t$, optional w.r.t. the
filtration $\widehat {\mathbb F}$ on $\widehat \Omega=\Omega\times
\mathbb R^+$, generated by $\mathcal F_t\otimes \mathcal B(\mathbb
R^+)$, such that
\begin{itemize}
\item for any $\theta$, $\alpha_.(\theta)$ is an $\mathbb
F$-martingale; moreover, denoting $\zeta^\theta= \inf \{t :
\alpha_{t-}(\theta)=0 \}$, then $\alpha_.(\theta)>0$, and
$\alpha_-(\theta)>0$ on $[0,\zeta^\theta)$, and
$\alpha_.(\theta)=0$ on $[\zeta^\theta,\infty)$.
\item For any bounded   family $(Y_t(\omega,\theta),t\geq 0)$
measurable w.r.t. $\mathcal P(\mathbb F) \otimes \mathcal
B(\mathbb R^+)$, (where $\mathcal P(\mathbb F)$ is the $\mathbb
F$-predictable $\sigma$-field), the $\mathbb F$-predictable
projection of  the process $Y_t(\omega,\tau)$ is the process
$Y_t^{(p)}=\int \alpha_{t-}(\theta)Y_t(\theta)\eta(d\theta).$\\ In
particular, for any $t$,
$\proba(\zeta^\tau<t)=\esp[\int_0^\infty\alpha_{t-}(\theta)\indic_{\{\zeta^\theta<t\}}\eta(d\theta)]=0.$
So, $\zeta^\tau$ is infinite a.s.
\end{itemize}%
\vspace*{0,5cm}

  We  are also concerned with the c\`adl\`ag version of the
martingale $(S_t(u),t\geq 0)$ for any $u\in \mathbb R^+$. By the
previous result, we have a universal version of their predictable
projections,
$$S_{t-}(u)=S_t^{(p)}(u)=\int_u^\infty
\alpha_{t-}(\theta)\eta(d\theta).$$
It remains to define  $S_t(u)=\displaystyle\lim_{q\in {\mathbb
Q}^+,\, q\downarrow t}S_q^{(p)}(u)$ to obtain a
universal c\`adl\`ag version of the martingales $S_.(u)$.\\
Remark that to show directly that $\int_u^\infty
\alpha_{t}(\theta)\eta(d\theta)$ is a c\`adl\`ag process, we need
stronger assumption on the process $\alpha_{t}(\theta)$ which
allows us to apply the Lebesgue theorem w.r.t. $\eta(d\theta)$.

 We say that the process  $(Y_{t}(\omega,\theta),t\geq 0)$ is
$\mathbb F$-optional if it is $\mathcal O(\mathbb
F)\otimes\mathcal B(\R^+)$-measurable, where $\mathcal O(\mathbb
F)$ is the optional $\sigma$-field of $\mathbb F$. In particular,
the process $(Y_t(\omega,t),t\geq 0)$ is optional.

\subsection{Density and intensity processes}

We are now interested in the relationship between the density and
the  intensity process  of $\tau$.  As we shall see, this is
closely related to the (additive and multiplicative)
decompositions of the supermartingale $S$.
%
\subsubsection{$\mathbb F$-decompositions of the survival process
$S$}
%
 In this section, we characterize the martingale and the predictable increasing part  of
the additive and multiplicative Doob-Meyer decomposition of the
supermartingale $S$ in terms of the density.
\begin{Pro}\label{Pro:decomposition S} {~}\\
1) The Doob-Meyer decomposition of  the survival process $S$  is
given by $S_t=1+M_t^{\mathbb F}-\int_0^t\alpha_u(u)\eta(du)$ where
$M^{\mathbb F}$ is the c\`adl\`ag  square-integrable
$\ff$-martingale defined by
\[\,M_t^{\mathbb F} =-\int_0^t\big(\alpha_t(u)-\alpha_u(u)\big)\eta(du)=
\esp[\int_0^\infty\alpha_u(u)\eta(du)|\F_t]-1,\; a.s.\]
2) Let $\zeta^{\mathbb F}:=\inf\{t: S_{t-}=0\}$ and define
$\lambda_{t}^{\mathbb F}:=\frac{\alpha_t(t)}{S_{t-}}$ for any $t <
\zeta^{\mathbb F}$ and  $\lambda_t^{\mathbb F}:=\lambda^{\mathbb
F}_{t\wedge\zeta^{\mathbb F}}$ for any $t\geq\zeta^{\mathbb F}$.
The multiplicative decomposition of $S$ is given by
\begin{equation}\label{St}
S_t=L_t^{\mathbb F}e^{-\int_0^t\lambda_s^{\mathbb F}\eta(ds)}\end{equation}
where  $L^{\mathbb F}$ is the $\mathbb F$-local martingale
solution of  $\hspace{3mm} dL_t^{\mathbb
F}=e^{\int_0^t\lambda_s^{\mathbb F}\eta(ds)}dM_t^{\mathbb F},\>
L_0^{\mathbb F}=1.$
\end{Pro}

\proof 1)  First notice that $(\int_0^t\alpha_u(u)\eta(du),t\geq
0)$ is an $\mathbb F$-adapted continuous increasing process (the
measure $\eta$ does not have any atom). By the martingale property
of $(\alpha_t(\theta),t\geq 0)$, for any fixed $t$,
 $$S_t=\proba (\tau >t\vert \F_t)=\int_t^\infty\alpha_t(u)\eta(du)
 =\esp[\int_t^\infty\alpha_u(u)\eta(du)|\F_t],\,\, a.s..$$
 Therefore,
the   non-negative process
$S_t+\int_0^t\alpha_u(u)\eta(du)=\esp[\int_
0^\infty\alpha_u(u)\eta(du)|\F_t]$ is a square-integrable
martingale since
$$\esp\Big[ \big(\int_0^\infty\alpha_u(u)\eta(du)\big)^2\Big]=2\esp\Big[\int_0^\infty\int_u^\infty
\alpha_s(s)\eta(ds)\alpha_u(u)\eta(du)\Big]=2\esp\Big[\int_0^\infty
S_u \alpha_u(u)\eta(du)\Big]\leq 2.$$  We shall choose its
c\`adl\`ag version if needed. Using the fact that $\int_0^\infty
\alpha_t(u)\eta (du)=1$, we obtain
$$\forall t,\;\esp[\int_ 0^\infty\alpha_u(u)\eta(du)|\F_t]=
1-\int_0^t(\alpha_t(\theta)-\alpha_\theta(\theta))\eta(d\theta),\;
a.s.$$
and the result follows.\\
2) Setting $L_t^{\mathbb F}=S_t e^{\int_0^t\lambda_s^{\mathbb
F}\eta(ds)}$, integration by parts formula  and 1) yield to
$$dL_t^{\mathbb F}=e^{\int_0^t\lambda_s^{\mathbb
F}\eta(ds)}dS_t+e^{\int_0^t\lambda_s^{\mathbb
F}\eta(ds)}\lambda_t^{\mathbb
F}S_t\eta(dt)=e^{\int_0^t\lambda_s^{\mathbb
F}\eta(ds)}dM_t^{\mathbb F},$$ which implies the result.
 \finproof

\begin{Rems}\label{remaimmer}{~}\\1) Note that, from (\ref{setA}),
${\mathbb P}(\zeta^{\mathbb F} \geq \tau)=1$. \\2) The survival
process $S$ is a decreasing process if and only  if the martingale
$M^{\mathbb F}$ is constant ($M^{\mathbb F}\equiv 0$) or
equivalently if and only if the martingale $L^{\mathbb F}$ is
constant ($L^{\mathbb F}\equiv 1$).  In that case, by Proposition
\ref{Pro:decomposition S}, $S$ is the continuous decreasing
process $S_t=e^{-\int_0^t\lambda_s^{\mathbb F}\eta(ds)}$.
Moreover, for any pair $(t,\theta),\, t\leq \theta$, the conditional
distribution is given by
$S_t(\theta)=\esp[e^{-\int_0^\theta\lambda_s^{\mathbb
F}\eta(ds)}|\F_t]$.\\
 3) The condition $M^{\mathbb F} \equiv 0$ can be written as
$\int_0^t (\alpha_t(u)-\alpha_u(u))\eta(du)=0$ and is satisfied
if, for $t\geq u$, $\alpha_t(u)-\alpha_u(u) =0$ (immersion
property), but the converse is obviously not true (Remark
\ref{pseudo}).
\end{Rems}

\subsubsection{Relationship with the $\mathbb G$-intensity}
\label{Subsec:density and intensity}

The intensity approach has been largely studied in the credit
literature. We study now in more details the relationship between
the density and the intensity, and notably between the $\mathbb
F$-density process of $\tau$ and its intensity process with
respect to $\mathbb G$. We first recall some definitions.
\begin{Def}
Let $\tau$ be a $\mathbb G$-stopping time. The {\it $\mathbb
G$-compensator}  of $\tau$ is the $\mathbb
G$-predictable increasing process $\Lambda^{\mathbb G}$ such that the process $(
N^{\mathbb G}_t=\indic_{\{\tau\leq t\}}-\Lambda^{\mathbb
G}_t,t\geq 0)$ is a $\mathbb G$-martingale. If $\Lambda^{\mathbb
G}$ is absolutely continuous with respect to the measure $\eta$,
the $\mathbb G$-adapted process $\lambda^{\mathbb G}$ such that
$\Lambda^{\mathbb G}_t=\int_0^t\lambda^{\mathbb G}_s\eta(ds)$ is
called the {\it $(\mathbb G,\eta)$-intensity process}  or the {\it
{$\mathbb G$-intensity}} if there is no ambiguity. The $\mathbb
G$-compensator is stopped at $\tau$, i.e., $\Lambda^{\mathbb
G}_t=\Lambda^{\mathbb G}_{t\wedge\tau}$.  Hence,
$\lambda_t^{\mathbb G}=0$ when $t>\tau$.
\end{Def}
\noindent The following results give the $\mathbb G$-intensity of
$\tau$ in terms of $\mathbb F$-density, and conversely the
$\mathbb F$-density $\alpha_t(\theta)$ in terms of the $\mathbb
G$-intensity, but only for $\theta\geq t$.
\begin{Pro}\label{Pro:density and intensity}
 {~}\\
1) The random time $\tau$ admits a $(\mathbb G, \eta)$-intensity
given by
\begin{equation}\label{Equ:compensator general case with F and G}\lambda_t^{\mathbb
G}=\indic_{\{\tau>t\}}\lambda_t^{\mathbb
F}=\indic_{\{\tau>t\}}\frac{\alpha_t(t)}{S_{t}},\quad \eta(dt)\>
a.s. .\end{equation}   The processes $(N^{\mathbb
G}_t:=\indic_{\{\tau\leq t\}}-\int_0^t\lambda^{\mathbb
G}_s\eta(ds),t\geq 0)$, and $(L^{\mathbb G}_t:=\indic_{\{\tau>
t\}}e^{\int_0^t\lambda_s^{\mathbb G}\eta(ds)},t\geq 0)$ are $\mathbb G$-local martingales.   \\
2) For any $t<\zeta^{\mathbb F}$ and $\theta \geq t$, we have:
$\alpha_t(\theta)=\esp[\lambda_\theta^{\mathbb G}|\F_t]$.\\
 Then, the $\mathbb F$-optional projections of the local
martingales $N^{\mathbb G}$ and $L^{\mathbb G}$ are the $\mathbb
F$-local martingales $-M^{\mathbb F}$ and $ L^{\mathbb F}$.
\end{Pro}
 \proof 1) The $\mathbb G$-local martingale property of
$  N^{\mathbb G}$   is equivalent to the $\mathbb G$-local
martingale property of $L^{\mathbb G}_t=\indic_{\{\tau>
t\}}e^{\int_0^t\lambda_s^{\mathbb G}\eta(ds)}=\indic_{\{\tau>
t\}}e^{\int_0^t\lambda_s^{\mathbb F}\eta(ds)}$, since
$$dL^{\mathbb G}_t=-L^{\mathbb G}_{t-}d\ind_{\{\tau\leq
t\}}+\ind_{\{\tau>t\}}e^{\int_0^t\lambda_s^{\mathbb
F}\eta(ds)}\lambda_t^{\mathbb F}\eta(dt)=-L^{\mathbb
G}_{t-}dN^{\mathbb G}_t$$
 Since the process $\int_0^t\lambda_s^{\mathbb F}\eta(ds)$ is
 continuous, we can proceed by localization introducing the $\mathbb G$-stopping times $\tau_n= \tau
\indic_{\{\tau \leq T_n\}} +\infty \indic_{\{\tau
> T_n\}}$ where $T_n=\inf\{t: \int_0^t\lambda_s^{\mathbb F}\eta(ds)
>n\}$.
Then, the martingale property of the stopped process $L^{\mathbb
G}_{t\wedge \tau^n}=L^{\mathbb G,n}_t$ follows from the $\mathbb
F$-martingale  property of
$L^{\mathbb F}_{t\wedge T^n}=L^{\mathbb F,n}_t$, since for any
$s\leq t$,
\begin{eqnarray*}
 && \esp[L^{\mathbb G,n}_t|\G_s]=\esp[\indic_{\{\tau> t\wedge T^n\}}e^{\int_0^{t\wedge
T^n}\lambda_u^{\mathbb F}\eta(du)}|\G_s]= \indic_{\{\tau> s\wedge
T^n\}}\frac{\esp[\indic_{\{\tau> t\wedge T^n\}}e^{\int_0^{t\wedge
T^n}\lambda_u^{\mathbb F}\eta(du)}|\F_s]}{S_s}\\
&&= \indic_{\{\tau> s\}}\frac{\esp[S_{t\wedge
T^n}e^{\int_0^{t\wedge T^n}\lambda_u^{\mathbb
F}\eta(du)}|\F_s]}{S_s}=\indic_{\{\tau> s\wedge
T^n\}}\frac{L_s^{\mathbb F,n}}{S_s}
\end{eqnarray*}
where the last equality follows from the $\mathbb F$-martingale
property of $L^{\mathbb F,n}$.\\ Then, the form of the intensities
follows from the definition.\\
2) By the martingale property of density, for any $\theta \geq t$,
$\alpha_t(\theta)=\esp[\alpha_\theta(\theta)|\F_t]$. Using the
definition of $S$, and the value of $\lambda ^{\mathbb G}$ given
in 1), we obtain
\[\alpha_t(\theta)=  \esp\Big[\alpha_\theta(\theta)
\frac{\indic_{\{\tau>\theta\}}}{S_{\theta}}|\F_t\Big]=\esp[\lambda_\theta^{\mathbb
G}|\F_t],\quad \forall t<\zeta^{\mathbb F},\; a.s.\]  Hence, the
value of the density can be partially deduced from the
intensity.\\
The $\mathbb F$-projection of the local martingale $L^{\mathbb
G}_t=\indic_{\{\tau> t\}}e^{\int_0^t\lambda_s^{\mathbb
G}\eta(ds)}$ is the local martingale
$S_te^{\int_0^t\lambda_s^{\mathbb F}\eta(ds)}=L^{\mathbb F}_t\>$
by definition of the survival process $S$. Similarly, since
$\alpha_t(\theta)=\esp[\lambda_\theta^{\mathbb G}|\F_t]$, the
$\mathbb F$-projection of the martingale $N^{\mathbb
G}_t=\indic_{\{\tau\leq t\}}-\int_0^t\lambda^{\mathbb
G}_s\eta(ds)$ is $1-S_t-\int_0^t\alpha_s(s)\eta(ds)=-M^{\mathbb
F}_t$. \finproof

\begin{Rems}$~$\\
1) Since the intensity process is continuous, $\tau$ is a totally
inaccessible $\mathbb G$-stopping time. \\
 2) The
density hypothesis, and the fact that $\eta$ is non-atomic allow
us to choose $\alpha_s(s)/S_s$ as an intensity, instead of
$\alpha_s(s)/S_{s-}$ as it is usually done (see \cite{EJY2000} in
the case where the numerator $\alpha_s(s)$ represents the
derivative of the compensator of $S$).  \\
3) Proposition \ref{Pro:decomposition S} shows that density and
intensity approaches correspond respectively to the additive and
the multiplicative decomposition point of view of the survival
process $S$.
\end{Rems}
We now use the density-intensity relationship to
  characterize
the pure jump $\mathbb G$-martingales   having only one jump at
$\tau$.

\begin{Cor}\label{Cor: G jump martingale} ~\\
1) For any  locally bounded $\mathbb G$-optional process
$H^{\mathbb G }$, the process
\begin{equation}\label{nh}N^{H,\mathbb G}_t:=H^{\mathbb G }_\tau \indic_{\{\tau \leq
t\}} -\int_0^{t\wedge \tau} \frac{\alpha_s(s)}{S_{s }}H^{\mathbb G
}_s\eta (ds)=\int_0^t\,H^{\mathbb G}_s\,dN^{\mathbb G}_s, \quad
t\geq 0
\end{equation}
 is a $\mathbb G$-local martingale.\\
  2) Conversely, any pure jump $\mathbb G$-martingale $M^{\mathbb G}$
which has only one locally bounded jump at $\tau$ can be written
on the form (\ref{nh}), with $H^{\mathbb G}_{\tau}=M^{\mathbb G}_{\tau}-M^{\mathbb G}_{\tau_-}$. \\
3)   Any nonnegative pure jump $\mathbb G$-martingale $U^{\mathbb
G}$ such that $U^{\mathbb G}_0=1$,  with only one jump at time
$\tau$ has the following representation
$$U^{\mathbb G}_t=\big(u_\tau\ind_{\{\tau\leq t\}}+\ind_{\{t<\tau\}}\big)
e^{-\int_0^{t\wedge \tau}(u_s-1)\lambda^{\mathbb F}_s\eta(ds)}$$
where $u $ is a positive $\mathbb F$-optional process associated
with the relative jump such that
 $u_{\tau}=U^{\mathbb G}_{\tau}/U^{\mathbb G}_{\tau-}$.
\end{Cor}
\proof  1) The $\mathbb G$-martingale property of $N^{\mathbb G}$
implies that $N^{H,\mathbb G}$ defined in \eqref{nh} is a $\mathbb
G$-martingale for any bounded predictable process $H^{\mathbb G}$.
From a reinforcement of (\ref{della}), if $H^{\mathbb G}$ is a
$\mathbb G$-predictable process (typically $H^{\mathbb G}_{t_0}
\indic_{]t_0,\infty]}$), there exists an $\mathbb F$-predictable
process $H^{\mathbb F}$ such that
$H_\tau^{\mathbb G} =H_\tau^{\mathbb F},\, a.s. .$
Then the process $H^{\mathbb G}$ may be   replaced by its
representative $H^{\mathbb F}$ in the previous relations.\\
Let $Y^{\mathbb G}_s$ be a bounded ${\mathcal G}_s$-random
variable, expressed in  terms of $\mathbb F$-random variables as $
Y^{\mathbb G}_s=Y^{\mathbb F}_s\ind_{\{s<\tau\}}+Y^{\mathbb
F}_s(\tau)\ind_{\{\tau\leq s\}}$ where $Y^{\mathbb F}_s\in
{\mathcal F}_s$ and $Y^{\mathbb F}_s(\theta)\in {\mathcal
F}_s\otimes {\mathcal B}([0,s])$, (typically $Y^{\mathbb
F}_s(\theta)=Y^{\mathbb F}_s\times g(\theta)\ind_{[0,s]})$. Then,
the ${\mathbb G}$-martingale property still holds for the process
$N^{H,{\mathbb G}}$ where $H^{\mathbb G}$ is the $\mathbb
G$-optional process $H^{\mathbb G}_s=Y^{\mathbb
G}_s\ind_{[s,\infty)}$. The optional $\sigma$-field being
generated by such processes, the assertion   holds for any
$\mathbb G$-optional process.\\
 2)  For the  converse, observe that
 the locally bounded jump $H^{\mathbb G}_{\tau}$ of the martingale $M^{\mathbb G}$ at
time $\tau$ is the value at time $\tau$ of some locally bounded
$\mathbb F$-optional process $H^{\mathbb F}$. Then the difference
$M^{\mathbb G}-N^{H,{\mathbb G}}$ is a finite variation local
martingale without jump, that is a constant process. \\
3) It is easy to calculate the differential of the finite
variation process $U^{\mathbb G} $ as
$$dU^{\mathbb G}_t=-U^{\mathbb G}_t(u_t-1)\lambda^{\mathbb G}_t\eta(dt)+U^{\mathbb G}_{t-}((u_t-1)(dN^{\mathbb
G}_t+\lambda^{\mathbb G}_t\eta(dt))=U^{\mathbb
G}_{t-}(u_t-1)dN^{\mathbb G}_t.$$ Then $U^{\mathbb G} $ is the
exponential martingale of the purely jump martingale
$(u_t-1)dN^{\mathbb G}_t$.
 \finproof
\subsection{An example of HJM type}
 We now give some examples, where we point out similarities with
 Heath-Jarrow-Morton models. Here, our
 aim is not  to  present a
 general framework, therefore, we reduce our attention to the
 case where the reference filtration $\mathbb F$ is generated by a  multidimensional  standard Brownian
 motion $W$.
The following two propositions, which model the dynamics of the
conditional
 probability  $S (\theta)$, correspond respectively to the additive and multiplicative
 points of view.
\noindent  From the predictable representation  theorem
 in the Brownian filtration, applied to the family
of bounded martingales $(S_t(\theta),t\geq 0)$, $\theta\geq 0$,
there exists a family of $\mathbb F$-predictable processes
$(Z_t(\theta),t\geq 0)$ such that
\begin{equation}\label{eq:survivalmartrepres}
dS_t(\theta)=Z_t(\theta)dW_t,\quad a.s.
\end{equation}
\begin{Pro}\label{Pro:HJM additive}  Let $dS_t(\theta)=Z_t(\theta)dW_t$
be the martingale representation of $(S_t(\theta),t\geq 0)$ and
assume that the processes $(Z _t(\theta);t\geq 0)$ are
differentiable in the following sense: there exists a family of
 processes
$(z_t(\theta),t\geq 0)$, bounded by an integrable process,  such
that $Z_t(\theta)=\int_0^{\theta}z_t(u)\eta(du)$. Then,

1) The density martingales have the following dynamics $\>
d\alpha_t(\theta)=-z_t(\theta)dW_t$.

2) The survival process $S$ evolves as
 $\> dS_t=-\alpha_t(t) \eta(dt) +\,Z_t(t)dW_t$.

 3) With more regularity assumptions,  if $(\partial_\theta
\alpha_t(\theta))_{\theta=t}$ is simply denoted by
 $\partial_\theta \alpha_t(t)$,  then the process $\alpha_t(t)$  is driven by :
 $$\> d\alpha_t(t)=\partial_\theta \alpha_t(t) \eta(dt)-z_t(t)dW_t$$
 \end{Pro}
\proof Observe  that   $Z (0)= 0$ since $S (0)= 1$, hence the
existence of $z$ is related with some smoothness conditions. Then
$$S_t(\theta)= S_0(\theta)+\int_0^t Z_u(\theta)dW_u= S_0(\theta) +
\int_0^\theta \eta(dv)\int_0^t z_u(v)dW_u$$ and 1) follows.
Furthermore,   by using Proposition \ref{Pro:decomposition S} and
integration by parts,
\[M_t^{\mathbb
F}=\int_0^t(\alpha_t(u)-\alpha_u(u))\eta(du)=\int_0^t\eta(du)\int_u^t
z_s(u)dW_s=\int_0^tdW_s\Big(\int_0^sz_s(u)\eta(du)\Big)
\] which implies 2).\\
3)   Let us use the short notation introduced above. We
follow the same way as for the decomposition of $S$, by studying
the process
$$\alpha_t(t)-\int_0^t (\partial_\theta \alpha_s)(s)\eta(d\theta )=\alpha_t(0)+
\int_0^t (\partial_\theta \alpha_t)(s)\eta(ds )-\int_0^t
(\partial_\theta \alpha_s)(s)\eta(ds)$$
Using martingale representation of $\alpha_t(\theta)$ and
integration by parts, (assuming that smoothness hypothesis allows
these operations) the integral in the RHS is a stochastic
integral,
\begin{eqnarray*}
&&\int_0^t \Big((\partial_\theta \alpha_t)(s)-(\partial_\theta
\alpha_s)(s)\Big)\eta(ds)=-\int_0^t
\eta(ds)\partial_\theta(\int_s^t\>z_u(\theta)dW_u)\\
&&=-\int_0^t \eta(ds)\int_s^t\>\partial_\theta
z_u(s)dW_u=-\int_0^t dW_u\int_0^u\eta(ds)\partial_\theta
z_u(s)=-\int_0^t dW_u (z_u(u)-z_u(0))
\end{eqnarray*}
The stochastic integral $\int_0^t dW_u z_u(0)$ is the stochastic
part of the martingale $\alpha_t(0)$, and so the property 3) holds
true.
 \finproof
We now consider $(S_t(\theta),t\geq 0)$ in the classical HJM
models (see \cite{Sch}) where its dynamics is given in
multiplicative form. By definition, the  forward hazard rate
$\lambda_t(\theta)$ of $\tau$ is given by
$\lambda_t(\theta)=-\partial_{\theta}\ln S_t(\theta)$ and the
density can then be calculated as
$\alpha_t(\theta)=\lambda_t(\theta)S_t(\theta)$. As noted in
Remark \ref{forward rate}, $\lambda(\theta)$ plays the same role
as the spot forward rate in the interest rate models.
\\
Classically, HJM framework is studied for time smaller than maturity, i.e. $t\leq T$. Here we
consider all positive pairs $(t,\theta)$.
\begin{Pro}\label{Pro:HJM multiplicative}
For any $t,\theta\geq 0$, let
$\Psi_t(\theta)=\frac{Z_t(\theta)}{S_t(\theta)}$ with the notation
of Proposition \ref{Pro:HJM additive}. We assume that
$\psi_t(\theta)$ defined by
$\Psi_t(\theta)=\int_0^{\theta}\psi_t(u)\eta(du)$ is bounded by
some integrable process. Then\\[2mm] 1)
$S_t(\theta)=S_0(\theta)\exp\left(\int_0^t\Psi_s(\theta)
dW_s-\frac12\int_0^t|\Psi_s(\theta)|^2ds\right)$;\\[2mm]
 2) The forward hazard  rate $\lambda(\theta)$ has the dynamics:
$\lambda_t(\theta)=\lambda_0(\theta)-\int_0^t \psi_s(\theta)dW_s
+\int_0^t\psi_s(\theta)\Psi_s(\theta)^*ds$;
 \\[2mm]
3) $S_t=\exp\left(-\int_0^t\lambda_s^{\mathbb F}\eta(ds)+\int_0^t
\Psi_s(s)dW_s-\frac{1}{2}\int_0^t|\Psi_s(s)|^2ds\right)$; 
 \end{Pro}
\proof By choice of notation, the process $S_t(\theta)$ is the
solution of the equation
\begin{equation}\label{Equ:equation of MtT}
\frac{dS_t(\theta)}{S_t(\theta)}=\Psi_t(\theta)dW_t,\qquad
\forall\, t,\theta\geq 0.
\end{equation}
Hence 1), from which we deduce immediately 2) by differentiation
w.r.t. $\theta$.\\
3) This representation is the multiplicative version of the
additive decomposition of $S$. There is not technical difficulties
because $S$ is continuous. \finproof
 \begin{Rems}\label{immer} If
$\Psi_s(s)=0$, then $S_t=\exp(-\int_0^t\lambda_s^{\mathbb
F}\eta(ds))$, which is decreasing. For the (H)-hypothesis to hold,
it needs $\Psi_s(\theta)=0$ for any $s\geq\theta$. \\ As a
conditional survival probability, $S_t(\theta)$ is decreasing on
$\theta$, which is equivalent to that $\lambda_t(\theta)$ is
positive.   When $\theta>t$, this property is implied by the
weaker condition $\lambda_t(t)\geq 0$. That is similar as for the
zero coupon bond prices. But when $\theta<t$, additional
assumption is necessary. We do not characterize this
condition.\end{Rems}
\begin{Rem}The above results are not restricted to the
Brownian filtration and can be easily extended to more general
filtrations under similar representation
$dS_t(\theta)=Z_t(\theta)dM_t$ where $M$ is a martingale which can
include  jumps. In this case, Proposition \ref{Pro:HJM additive}
can be generalized in a similar form; for Proposition \ref{Pro:HJM
multiplicative}, more attention should be payed to Dol\' eans-Dade
exponential martingales with jumps.
\end{Rem}
{\bf Example:}
We now give a particular example which provides a large class of
forward rate processes. The
non-negativity of $\lambda$ is satisfied, by 2) of Proposition \ref{Pro:HJM multiplicative}, if \\
\hspace*{1cm} $\bullet$ for any $\theta$, the process
$\psi(\theta)\Psi( \theta) $ is non
negative, or if $\psi( \theta)$ is non negative; \\
\hspace*{1cm} $\bullet$ for any $\theta$,  the
 local  martingale
$\zeta_t(\theta)=\lambda_0(\theta)-\int _0^t \psi_s(\theta)  dW_s$
is a Dol\'eans-Dade exponential of some martingale, i.e., is
 solution
of $$ \zeta_t (\theta) = \lambda_0(\theta)+\int_0^t \zeta_s
(\theta)  b_s(\theta) dW_s\, ,$$ that is, if $-\int _0^t
\psi_s(\theta) dW_s= \int_ 0^t b_s( \theta) \zeta_s (\theta) d
W_s$. Here the initial condition is a positive constant
$\lambda_0(\theta) $. Hence, we set
$$\psi _t (\theta)=  - b_t (\theta) \zeta_t (\theta)= -b_t (\theta)\lambda_0(\theta)\exp
\left ( \int_0^t b_s(\theta) d W_s-\frac 12\int_0^t  b_s^2(\theta)
ds\right)$$  where $ \lambda_0 $ is  a positive intensity function
and $b (\theta)$ is  a non-positive  $\ff$-adapted process. Then,
the family  $$ \alpha_t(\theta)=\lambda_t(\theta)\exp \left(-
\int_0^{\theta } \lambda_t(v)\, dv\right), $$ where
\begin{eqnarray*}
\lambda_t(\theta) =\lambda_0(\theta)- \int_0^t \psi_s(\theta)\,
dW_s+\int_0^t\psi_s(\theta)\Psi_s (\theta) \, ds
\end{eqnarray*}
satisfies the required assumptions.
\section{Characterization of $\mathbb G$-martingales in terms of $\mathbb F$-martingales}

In the theory of pricing and hedging,  martingale properties play
a very important role. In this section, we study the martingale
characterization when taking into account information of  the
default occurrence. The classical question in the enlargement of
filtration  theory is to give decomposition of $\mathbb
F$-martingales in terms of $\mathbb G$-semimartingales. For the
credit problems, we are concerned with the problem in a converse
sense, that is, with the links between   $\mathbb G$-martingales
and $\mathbb F$-(local) martingales.  In the literature, $\mathbb
G$-martingales which are stopped at $\tau$ have   been
investigated, particularly in the credit context. For our analysis
of after-default events, we are furthermore interested in the
martingales which start at the default time $\tau$ and in
martingales having one jump at $\tau$, as the ones introduced in
Corollary \ref{Cor: G jump martingale}.  We shall give
characterization results for these types of $\mathbb
G$-martingales in the following, by using a coherent formulation
in the density framework.

\subsection{$\mathbb G$-martingale characterization}
\label{Subsec:martingale representation one name}
%
Any $\mathbb G$-martingale may be split into two martingales, the
first one stopped at time $\tau$ and the second one starting at
time $\tau$, that is \[Y_t^{\mathbb G}=Y_t^{\rm{bd},\mathbb
G}+Y_t^{\rm{ad},\mathbb G}\] where $Y_t^{\rm{bd},\mathbb
G}:=Y^{\mathbb G}_{t\wedge\tau}$ and $Y_t^{\rm{ad},\mathbb G}:=(
Y_t^{\mathbb G}-Y_\tau^{\mathbb G})\indic_{\{\tau\leq t\}}$. We
now study the two types of martingales respectively.
\\
The density hypothesis allows us to provide easily a
characterization\footnote{The following proposition was
established in \cite[Lemma 4.1.3]{BJRanta} in a hazard process
setting.} of $\mathbb G$-martingales stopped at time $\tau$.

\begin{Pro}\label{Pro:martingale stopped at tau}
A $\mathbb G$-adapted c\`adl\`ag process $Y^{\mathbb G}$ is a
closed $\mathbb G$-martingale stopped at time $\tau$ if and only
if there exist an $\mathbb F$-adapted c\`adl\`ag process $Y$
defined on $[0,\zeta^{\mathbb F})$ and an $\mathbb F$-optional
process $Z$ such that $Y^{\mathbb G}_t=Y_t\indic_{\{\tau>t\}}+
Z_{\tau} \indic_{\{\tau\leq t\}}$ a.s. and that
\begin{equation}\label{Equ:cond mart stopped at tau}
(U_t:=Y_tS_t+\int_0^tZ_s \alpha_s(s)\eta(ds),\,t\geq 0)
\quad\text{is an} \,\,\mathbb F\text{-martingale on}\;
[0,\zeta^{\mathbb F}).
\end{equation}
Equivalently, using the multiplicative decomposition of $S$ as
$S_t=L_t^{\mathbb F}e^{-\int_0^t\lambda_s^{\mathbb F}\eta(ds)}$ on
$[0,\zeta^{\mathbb F})$, the above condition \eqref{Equ:cond mart
stopped at tau} is equivalent to
\begin{equation}\label{Equ:second cond mart stopped at tau}
(L_t^{\mathbb F}[Y_t+\int_0^t(Z_s -Y_s)\lambda_s^{\mathbb
F}\eta(ds)],\,t\geq 0) \,\,\text{is an} \,\,\mathbb F\text{-local
martingale on}\; [0,\zeta^{\mathbb F}).
\end{equation}
\end{Pro}
\proof The conditional expectation of $Y_t^{\mathbb G}$ given
$\F_t$ is the $\mathbb F$-martingale defined on $[0,\zeta^{\mathbb
F})$ as $Y_t^{\mathbb F}=\esp[Y_t^{\mathbb
G}|\F_t]=Y_tS_t+\int_0^t Z_s \alpha_t(s)\eta(ds)$ by using the
$\F_t$-density of $\tau$. Notice that $Y_t^{\mathbb F}$ differs
from $U_t$ by $(\int_0^t Z_s
(\alpha_t(s)-\alpha_s(s))\eta(ds),t\geq 0)$, which is an
   $\mathbb F$-local martingale (this can be easily checked using
   that
$Z$ is locally bounded and $(\alpha_t(s),t\geq 0)$ is $\mathbb
F$-martingale). So $U $ is also an $\mathbb F$-local martingale.
Moreover, since $\E[ \vert Y_t^{\mathbb G} \vert ]<\infty$, for
any $\mathbb F$-stopping time $\vartheta$, the quantity
$Y_\vartheta\indic_{\{\tau>\vartheta\}}$ is integrable, hence
$Y_\vartheta S_\vartheta$ is also integrable, and $$\E[ \int_
0^{\zeta ^{\mathbb F}} \vert Z_s\vert \alpha_s(s)\eta (ds)]= \E
[\vert Y_\tau ^{\mathbb G}\vert]<\infty,$$ which establishes that
$U$ is
a martingale.\\
Conversely, if $U$ is an $\mathbb F$-local martingale, it is easy
to verify by Theorem \ref{Thm:cond expect general case with
density} that $\esp[Y_T^{\mathbb G}-Y_t^{\mathbb G}|\G_t]=0$,
a.s..\\
The second formulation is based on the multiplicative
representation $S_t=L_t^{\mathbb F}e^{-\Lambda_t^{\mathbb F}}$
where $\Lambda_t^{\mathbb F}=\int_0^t\lambda_s^{\mathbb
F}\eta(ds)$ is a continuous increasing process. Since
$e^{\Lambda_t^{\mathbb F}}Y_tS_t=Y_tL_t^{\mathbb F}$ and
$\alpha_t(t)=\lambda_t^{\mathbb F}S_{t}$, we have
\[d(Y_tL_t^{\mathbb F})=e^{\Lambda_t^{\mathbb F}}d(Y_tS_t)+e^{\Lambda_t^{\mathbb F}}Y_tS_t\lambda_t^{\mathbb
F}\eta(dt) =e^{\Lambda_t^{\mathbb F}}dU_t+(Y_t-Z_t
)\lambda_t^{\mathbb F}L_t^{\mathbb F}\eta(dt).\] The local
martingale property of the process $U$ is then equivalent to that
of $(Y_tL_t^{\mathbb F}-\int_0^t(Y_s- Z_s )\lambda_s^{\mathbb
F}L_s^{\mathbb F}\eta(ds),t\geq 0)$, and then to the condition
\eqref{Equ:second cond mart stopped at tau}.  \finproof
\begin{Rem}\label{Rem:continuite on tau}A $\mathbb G$-martingale
stopped at time $\tau$ and equal to 1 on $[0,\tau)$ is constant on
$[0,\tau]$. Indeed, integration by parts formula  proves that
$(L_t^{\mathbb F}\int_0^t(1-Z_s )\lambda_s^{\mathbb
F}\eta(ds),t\geq 0)$ is a local martingale if and only if   the
continuous bounded variation process $( \int_0^t L_s^{\mathbb F}
(1-Z_s )\lambda_s^{\mathbb F}\eta(ds),t\geq 0)$ is a
local-martingale, that is if $ L_s^{\mathbb F} (1-Z_s
)\lambda_s^{\mathbb F} =0$, which implies that $Z_s =1$ on
$[0,\zeta^{\mathbb F})$.
\end{Rem}
\noindent The before-default $\mathbb G$-martingale
$Y^{\rm{bd},\mathbb G}$ can always be separated into two parts: a
martingale which is stopped at $\tau$ and is continuous at $\tau$;
and a martingale which has a jump at $\tau$.
\begin{Lem} \label{Pro:martingale continue et saut pur}
Let $Y^{\rm{bd},\mathbb G}$ be a $\mathbb G$-martingale stopped at
$\tau$ of the form $ Y^{\rm{bd},\mathbb
G}_t=Y_t\indic_{\{\tau>t\}}+ Z_{\tau} \indic_{\{\tau\leq t\}}.$
Then there exist two $\mathbb G$-martingales $Y^{\rm{c,bd}}$ and
$Y^{\rm{d,bd}}$ such that $Y^{\rm{bd},\mathbb
G}=Y^{\rm{c,bd}}+Y^{\rm{d,bd}}$ which satisfy the following
conditions:
\\
1) $(Y^{\rm{d,bd}}_t=(Z_{\tau} -Y_\tau)\indic_{\{\tau\leq
t\}}-\int_0^{t\wedge\tau}(Z_s -Y_s)\lambda_s^{\mathbb
F}\eta(ds),t\geq 0)$ is a   $\mathbb G$-martingale with  a single jump at $\tau$; \\
2) $(Y^{\rm{c,bd}}_t=\widetilde Y_{\tau\wedge t},t\geq 0)$ is
continuous at $\tau$, where $\widetilde
Y_t=Y_t+\int_0^{t}(Z_s-Y_s)\lambda_s^{\mathbb F}\eta(ds)$.
  \end{Lem}
\proof From Corollary \ref{Cor: G jump martingale},
$Y^{\rm{d,bd}}$ is a martingale. The result follows. \finproof

\begin{Cor}\label{Cor:martingale continue at tau} With the above notation, a
martingale $Y^{\mathbb G}$ which is stopped and continuous at
$\tau$ is characterized by: $(L^{\mathbb F}_tY_t,t\geq 0)$ is an
$\mathbb F$-local martingale. Furthermore, if $L^{\mathbb F}$ is a
martingale, then this condition is equivalent to that $Y$ is an
$\mathbb F$-local martingale w.r.t. the probability measure
$\proba^L=L_T^{\mathbb F}\proba$.  In particular, under the
immersion assumption, the $\mathbb G$-martingales stopped at time
$\tau$ and continuous at $\tau$ are $\mathbb F$-martingales stopped
at $\tau$.
\end{Cor}
\begin{Rem}Under immersion, $L_t^{\mathbb F}=1$.
So a process $Y^{\mathbb G}$ stopped at $\tau$ and continuous at
time $\tau$ is a $\mathbb G$-martingale if and only if $Y$ is an
$\mathbb F$-local martingale.
\end{Rem}

 We now concentrate on the $\mathbb G$-martingales
starting at $\tau$, which, as we can see below, are easier to
characterize. The following proposition is a direct consequence of
Theorem \ref{Thm:cond expect general case with density}.

\begin{Pro}\label{Pro:martingale starting at tau}
Any c\`adl\`ag  integrable process $Y^{\mathbb G}$ is a $\mathbb
G$-martingale starting at $\tau$ with $Y_\tau=0$ if and only if
there exists an $\mathcal O(\mathbb F)\otimes \mathcal
B(\R^+)$-optional process $(Y_t(.),t\geq 0)$ such that  $Y_t(t)=0$
and $Y_t^{\mathbb G}=Y_t(\tau)\indic_{\{\tau\leq t\}}$ and that,
for any $\theta
>0$, $(Y_t(\theta)\alpha_t(\theta),t\geq \theta\geq 0)$ are
$\mathbb F$-martingales on $[\theta,\zeta^{\theta})$, where
$\zeta^{\theta}$ is defined as in Section \ref{Subsec:path
regularity}.
\end{Pro}

\noindent Combining the previous   results, we give the
characterization of general $\mathbb G$-martingale.

\begin{Thm}\label{Thm:G maringale characterization}
A  c\`adl\`ag process $Y^{\mathbb G}$ is a $\mathbb G$-martingale
if and only if there exist an $\mathbb F$-adapted c\`adl\`ag
process $Y$ and an $\mathcal O(\mathbb F) \otimes \mathcal
B(\R^+)$-optional  process $Y_t(.)$ such that $Y_t^{\mathbb
G}=Y_t\indic_{\{\tau>t\}}+Y_t(\tau)\indic_{\{\tau\leq t\}}$ and
  \\1) the process
$(Y_tS_t+\int_0^tY_s(s)\alpha_s(s)\eta(ds),\,t\geq 0)$ or
equivalently $(L_t^{\mathbb
F}[Y_t+\int_0^t(Y_s(s)-Y_s)\lambda_s^{\mathbb F}\eta(ds)],t\geq
0)$  is an $\mathbb F$-local martingale; \\ 2) for any $\theta\geq
0$, $(Y_t(\theta)\alpha_t(\theta),t\geq \theta)$ is an $\mathbb
F$-martingale on $[\theta,\zeta^{\theta})$.
\end{Thm}

\proof Notice that $Y_t^{\rm{ad},\mathbb
G}=(Y_t(\tau)-Y_\tau(\tau))\indic_{\{\tau\leq t\}}$. Then the
theorem follows directly by applying Propositions
\ref{Pro:martingale stopped at tau} and \ref{Pro:martingale
starting at tau} on $Y^{\rm{bd},\mathbb G}$ and
$Y^{\rm{ad},\mathbb G}$ respectively.\finproof

\begin{Rem}{~}  We observe again the fact that to characterize what goes on
before the default, it suffices to know the survival process $S$
or the intensity $\lambda^{\mathbb F}$. However, for the
after-default studies, we need the whole conditional distribution
of $\tau$, i.e., $\alpha_t(\theta)$ where $\theta\leq t$.
\end{Rem}

\subsection{Decomposition of $\mathbb F$-(local) martingale}
An important result in the enlargement of filtration theory is
the decomposition of $\mathbb F$-(local) martingales as  $\mathbb
G$-semimartingales. Using the above results, we provide  an
alternative proof for a result established  in \cite{JLC2008},
simplified by using the fact that any $\mathbb F $-martingale is
continuous at time $\tau$. Our method is interesting, since it
gives the intuition of the decomposition without using any result
on enlargement of filtrations.

\begin{Pro}\label{Cor:F martingal decomposition}
Any $\ff$-martingale $Y^{\mathbb F}$ is a $\mathbb
G$-semimartingale can be written as  $Y^{\mathbb F} _{t } = M^{Y,
\mathbb G}_t+A ^{Y, \mathbb G}_t  $ where $M^{Y, \mathbb G}$ is a
$\mathbb G$-martingale  and $(A_t^{Y, \mathbb G}:= A_t
\indic_{\{\tau
>t\}}+ A_t(\tau) \indic_{\{\tau \leq t\}},t\geq 0)$ is an optional
process with finite variation. Here
\begin{equation}\label{G variation fini A}A_t= \int_0^{t } \frac{d[Y^{\mathbb F},S]_s}{S_ {s}}\,\,
\text{ and } \,\,A_t(\theta)= \int_\theta^t \frac{d[Y^{\mathbb
F},\alpha (\theta)]_s}{\alpha _{s}(\theta)}.\end{equation}
\end{Pro}
\proof On  the one hand,  assuming that $Y^{\mathbb F}$ is a
$\mathbb G$-semimartingale,   it can be decomposed as the sum of a
$\mathbb G$-(local)martingale and a $\mathbb G$-optional process
$A ^{Y, \mathbb G}$ with finite variation which can be written as
$A_t \indic_{\{\tau
>t\}}+ A_t(\tau) \indic_{\{\tau \leq t\}}$ where $A$ and $A (\theta)$ are still unknown.
Note that, since $Y^{\mathbb F}$ has no jump at $\tau$ (indeed,
$\tau$ avoids $\mathbb F$-stopping times - see Corollary
\ref{Cor:avoid F stopping time}), we can choose $M^{Y,\mathbb G}$
such that $M^{Y,\mathbb G}$ and hence $A^{Y,\mathbb G}$ have no
jump at $\tau$. Applying  the martingale characterization result
obtained in Theorem \ref{Thm:G maringale characterization} to the
$\mathbb G$-local martingale
$$Y^{\mathbb F} _{t }-A ^{Y, \mathbb G}_t=(Y^{\mathbb F} _{t }-
A_t) \indic_{\{\tau
>t\}}+ (Y^{\mathbb
F} _{t }-A_t(\tau)) \indic_{\{\tau \leq t\}}$$  leads to the fact
that the two processes \begin{equation}\label{semim}((Y_t^{\mathbb
F}-A_t)L^{\mathbb F}_t ,t\geq 0) \,\,\text{ and } \,\,(\alpha
_t(\theta) (Y_t^{\mathbb F} -A_t(\theta)),t\geq
\theta)\end{equation} are $\mathbb F$-(local) martingales. Since
$$d\big((Y_t^{\mathbb F}-A_t)L_t^{\mathbb F}\big)=
(Y_{t-}^{\mathbb F}-A_{t-})dL_t^{\mathbb F}+L_{t-}^{\mathbb F}
d(Y_t^{\mathbb F}-A_t) +d\langle Y^{\mathbb F},L^{\mathbb
F}\rangle ^c_t+\Delta(Y_t^{\mathbb F}-A_t)\Delta L_t^{\mathbb
F}$$and
$$-L_{t-}^{\mathbb F}dA_t-\Delta A_t\Delta L_t^{\mathbb F}=-L_t^{\mathbb F}dA_t,$$
based on the intuition given by the Girsanov theorem, natural
candidate for the finite variation processes $A$ is
$dA_t=d[Y^{\mathbb F},L^{\mathbb F}]_t/L^{\mathbb F}_t$ where $[
~,~]$ denotes the co-variation process. Similarly,
$dA_t(\theta)=d[Y^{\mathbb
F},\alpha(\theta)]_t/\alpha_{t}(\theta)$.  Then,  using the fact
that   $Y^{\mathbb F}, L^{\mathbb F}, \alpha(\theta)$ are
${\mathbb F}$-local martingales, we obtain that $A=(1/L^{\mathbb
F}) \star [Y^{\mathbb F},L^{\mathbb F}]$ where $\star$ denotes the integration of $1/L^{\mathbb
F}$ w.r.t. $[Y^{\mathbb F},L^{\mathbb F}]$,  and
$A(\theta)=(1/\alpha(\theta)) \star [Y^{\mathbb
F},\alpha(\theta)]$ similarly. Then, since  $S$ is the product of the
martingale $L^F$ and an continuous increasing process
$e^{\Lambda^{\mathbb F}}$,
we have $d[Y^{\mathbb F},L^{\mathbb F}]_t/L_t^{\mathbb F}=d[Y^{\mathbb F},S]_t/S_t$ and obtain the first equality in (\ref{G variation fini A}).  \\
On the other hand, define the optional  process $A^{Y, \mathbb G}$
by using \eqref{G variation fini A}. It is not difficult to verify
by Theorem \ref{Thm:G maringale characterization} that $Y^{\mathbb
F} -A ^{Y, \mathbb G} $ is a $\mathbb G$-local martingale. It
follows that $Y^{\mathbb F}$ is indeed a $\mathbb
G$-semimartingale. \finproof
\begin{Rem} Note that our decomposition differs from the usual
one, since our process $A$ is optional (and not predictable) and
that we are using the co-variation process, instead of the
predictable co-variation process. As a consequence our
decomposition is not unique.\end{Rem}

\subsection{Girsanov theorem}

Change of probability measure is a key tool in derivative pricing
as in martingale theory. In credit risk framework, we are also
able to calculate parameters of the conditional distribution of
the default time w.r.t. a new probability measure. The links
between change of probability measure and the initial enlargement
have been established, in particular, in \cite{Ja1987} and
\cite{ADI2007}. In statistics, it is motivated by the Bayesian
approach \cite{GVV2006}.

 \noindent We present a Girsanov type
result, where the Radon-Nikod\'ym density is given in an additive
form instead of in a multiplicative one as in the classical
literature. This makes the density  of $\tau$ having simple form
under the new probability measure.

\begin{Thm}\label{Thm:Girsanov}\textbf{(Girsanov's theorem)}  Let $Q_t^{\mathbb
G}=q_t\indic_{\{\tau>t\}}+q_t(\tau)\indic_{\{\tau\leq t \}}$ be a
c\`adl\`ag positive $\mathbb G$-martingale with $Q_0^{\mathbb
G}=q_0=1$.   Let $\Q$ be the probability measure defined on $\G_t$
by $d\Q=Q_t^{\mathbb G}d\proba$ for any $t \in \R_+$ and  $\mathbb
Q^{\mathbb F}$ be the restriction of $\Q$ to $\mathbb F$, which
has Radon-Nikod\'ym density $Q^{\mathbb F}$, given by the
projection of $Q^{\mathbb G}$ on $\mathbb F$, that is
$Q_t^{\mathbb
F}=q_tS_t+\int_0^tq_t(u)\alpha_t(u)\eta(du)$.\\
Then $(\Omega,\Q,\mathbb G,\mathbb F,\tau)$ satisfies the density
hypothesis with the $(\mathbb F,\Q)$-density of $\tau$ given in
closed form after the default, that is for $\theta \leq t$ by
\[\alpha_t^{\Q}(\theta)=\alpha_t(\theta)
\frac{q_t(\theta)}{Q_t^{\mathbb F}}, \quad \;\eta(d\theta) \mbox{-
a.s.};\] and, only via a conditional expectation before the
default, that is for $t\leq \theta$, by
$$ \alpha_t^{\mathbb Q}(\theta)= \frac{1}{Q^{\mathbb F}_t}\,
\E^{\mathbb P}[ \alpha_\theta(\theta) q_\theta(\theta) \vert
\F_t].
$$
Furthermore: \\
1) the $\Q$-conditional survival process is defined on
$[0,\zeta^{\mathbb F})$ by
$S_t^{\Q}=S_t\displaystyle\frac{q_t}{Q_t^{\mathbb F}}$, and is
null after $\zeta^{\mathbb
F}$; \\
2) the
 $(\mathbb F,\Q )$-local martingale
$L^{\mathbb F,\Q}$ is   $(L_t^{\mathbb F,\Q}=L_t^{\mathbb
F}\displaystyle\frac{q_t}{Q_t^{\mathbb
F}}\exp{\int_0^t(\lambda_s^{\mathbb F,\Q}-\lambda_s^{\mathbb
F})\eta(ds)},t\in[0,\zeta^{\mathbb F}));$\\3) the
 $(\mathbb F,\Q)$-intensity process
is $\lambda_t^{\mathbb F,\Q}=\lambda_t^{\mathbb
F}\displaystyle\frac{q_t(t)}{q_{t}}, \eta(dt)$-a.s.
\end{Thm}
\proof The expression of the density process after the default
$(\alpha^{\Q}_t(\theta), \theta\leq t)$ is an immediate
consequence of definition. Before the default, the density may be
only obtained via a conditional expectation form given by
\begin{eqnarray*} \alpha_t^{\mathbb Q}(\theta)&=&
\E^{\mathbb Q}\big[\alpha_\theta(\theta)
\frac{q_\theta(\theta)}{Q_\theta^{\mathbb F}}\vert \F_t\big] =
\frac{1}{Q^{\mathbb F}_t}\;\E^{\mathbb P}\big[Q_\theta^{\mathbb F}
\alpha_\theta(\theta) \frac{q_\theta(\theta)}{Q_\theta^{\mathbb
F}}\vert \F_t\big]
\\&=&\frac{1}{Q^{\mathbb F}_t}\;\E^{\mathbb
P}[ \alpha_\theta(\theta)  q_\theta(\theta) \vert \F_t].
\end{eqnarray*}
For any $t\in[0,\zeta^{\mathbb F})$, the $\Q$-conditional survival
probability can be calculated
by\[S_t^{\Q}=\Q(\tau>t|\F_t)=\frac{\esp^{\proba}[\indic_{\{\tau>t\}}Q_t^{\mathbb
G}|\F_t]}{Q_t^{\mathbb F}}=q_t\frac{S_t}{Q_t^{\mathbb F}}\] and
finally, we use $\lambda_t^{\mathbb
F,\Q}=\alpha^{\Q}_t(t)/S^{\Q}_{t}$ and $L_t^{\mathbb F,\mathbb
Q}=S_t^{\Q}e^{\int_0^t\lambda_s^{\mathbb F,\Q}\eta(ds)}$ to
complete the proof. \finproof

It is known, from
\cite{GP2001}, that under density hypothesis, there exists at
least a change of probability, such that immersion property holds
under this change of probability.   Theorem \ref{Thm:Girsanov}
provides a full characterization of such changes of
probability.
\begin{Cor}\label{Cor:Changeto immersion}
We keep the notation of Theorem \ref{Thm:Girsanov}. The change of
probability measure generated by the two processes
\begin{equation}
q_t=(L^{\mathbb F}_t)^{-1},\qquad
q_t(\theta)=\frac{\alpha_\theta(\theta)}{\alpha_t(\theta)}
\end{equation} provides a model where the immersion property holds
true under $\Q$, and where the intensity processes  does not change, i.e., remains $\lambda^{\mathbb F}$.\\
More generally, the only changes of probability measure for which
the immersion property holds with the same intensity process are
generated by a  process $q$ such that $(q_t L^{\mathbb F}_t,t\geq 0)$ is
a uniformly integrable martingale.
\end{Cor}
\bproof: Any change of probability measure with immersion property
and the same intensity processes is characterized by the
martingale property of the product $Q^{\mathbb F}=q.\,L^{\mathbb
F}$. Moreover, given $q$, the immersion property determines in a
unique way the process $(q_t(\theta) ; t\geq \theta)$ via the
boundary condition $q_{\theta}(\theta)=q_{\theta}$ and the
equalities
$$
\alpha^{\mathbb
Q}_t(\theta)=\alpha_t(\theta)\frac{q_t(\theta)}{q_t\,L^{\mathbb
F}_t}=\alpha_{\theta}(\theta)\frac{q_{\theta}(\theta)}{q_{\theta}\,L^{\mathbb
F}_{\theta}} =\alpha_{\theta}(\theta)(1/L^{{\mathbb
F}}_{\theta}).$$
The martingale $Q^{\mathbb F}=q.\,L^{\mathbb F}$
has to satisfy the compatibility condition
\begin{eqnarray}
Q^{\mathbb F}_t&=&q_t\,L^{\mathbb F}_t=q_tS_t+Q^{\mathbb
F}_t \,\int_0^t\alpha_u(u)(1/L^{{\mathbb F}}_u)\eta(du) \nonumber\\
&=&Q^{\mathbb F}_t\big( e^{-\int_0^t\lambda^{\mathbb F,{\mathbb
P}}_s\eta(ds)} +\int_0^t e^{-\int_0^u\lambda^{\mathbb F, {\mathbb
P}}_s\eta(ds)} \lambda^{\mathbb F,{\mathbb P}}_u\>\eta(du)\big)
\label{compatibility} \end{eqnarray}
where the last equality comes from the identities \eqref{St} and \eqref{Equ:compensator general case with F and G}. The term in the bracket in \eqref{compatibility} is of finite variation and is hence equal to $1$. Then the $Q^{\mathbb F}$-compatibility condition is always satisfied. So
the only constraint on the process $q$ is the martingale property
of $q.\,L^{\mathbb F}$. \finproof

It is well known, from Kusuoka \cite{Ku1999}, that immersion
property is not stable by a change of probability. In the following, we shall in a first step characterize, under density hypothesis, changes of probability  which preserve
this immersion property, that is, H-hypothesis is satisfied under both $\proba$ and $\Q$. 
(See also \cite{CJN09} for a different study of changes of probabilities
preserving immersion property.) In a second step, we shall study
changes of probability which preserve the information before the
default, and give the impact of a change of probability after the
default.


\begin{Cor} \label{Cor:Girsanovimmersion}We keep the notation of Theorem \ref{Thm:Girsanov},
and assume immersion property under $\mathbb P$.\\
 1) Let the
Radon-Nikod\'ym density $(Q_t^{\mathbb G},t\geq 0)$ be a pure jump
martingale with only one jump at time $\tau$. Then, the $(\mathbb
F,\mathbb P)$-martingale $(Q_t^{\mathbb F},t\geq 0)$ is the
constant martingale equal to 1. Under $\Q$, the intensity process
is $\lambda_t^{\mathbb F,\Q}=\lambda_t^{\mathbb
F}\displaystyle\frac{q_t(t)}{q_{t}},
\eta(dt)$-a.s., and the immersion property still holds.\\
2) Conversely, the only changes of probability measure compatible
with  immersion property have  Radon-Nikod\'ym densities that are
the product of a pure jump positive martingale with only one jump
at time $\tau$, and a positive $\mathbb F$-martingale.
\end{Cor}

 \proof: The previous Girsanov  theorem \ref{Thm:Girsanov} gives immediately the intensity
characterization.\\
From Lemma \ref{Pro:martingale continue et saut pur}, the pure
jump martingale $(Q_t^{\mathbb G},t\geq 0)$ is a finite variation
process and   $(Q_t^{\mathbb G}=q_t, t<\tau)$ is a continuous
process with bounded variation. Since immersion property holds,
$S$ is a continuous decreasing process  (see footnote
\ref{scontinu}), and $Q_t^{\mathbb
F}=q_tS_t+\int_0^tq_u(u)\alpha_t(u)\eta(du)=q_tS_t+\int_0^tq_u(u)\alpha_u(u)\eta(du)
$ is a continuous martingale with finite variation. Since
$Q_0^{\mathbb F}=1$, then at any time, $Q_t^{\mathbb F}=1,\>a.s$.
By the other results established in Girsanov's theorem
\ref{Thm:Girsanov}, this key point implies that the new density is
constant after the default, so that the immersion property still
holds.\\
2) Thanks to the first part of this corollary, we can restrict our
attention to the case when in the both universe the intensity
processes are the same. Then the Radon-Nikod\'ym density is
continuous at time $\tau$ and the two processes $(q_t,t\geq 0)$
and $(q_t(\theta),t\geq \theta)$ are
$\mathbb F$-(local) martingales.\\
Assume now that the immersion property holds   also under the new
probability measure $\mathbb Q$. Both martingales $ L ^{\mathbb
F,\mathbb P}$ and $ L ^{\mathbb F,\mathbb Q}$ are constant, and
$Q^{\mathbb F}_t=q_t$. Moreover the $\mathbb Q$-density process
being constant after the default $(\theta<t)$,
$q_t(\theta)/q_t=q_\theta(\theta)/q_\theta=1, \>a.s.$. The
processes $Q^{\mathbb G}$, $Q^{\mathbb F}$ and $q$ are
undistinguishable. \finproof


 As shown in this paper, the knowledge of the intensity does
not allow to  give full information on the law of the default,
except if immersion property holds. Starting with a model under
which immersion property holds, taking $q_t(t)=q_t$ in Theorem
\ref{Thm:Girsanov} will lead us to a model where the default time
admits the same intensity whereas immersion property does not
hold, and then the impact of the default changes the dynamics of
the default-free assets. We present a specific case where, under
the two probability measures, the dynamics of these assets are the
same before the default but are changed after the default, a
phenomenon that is observed in the actual crisis.
%
 We   impose that the new probability $\mathbb Q$   coincide with
$\mathbb P$ on the $\sigma$-algebra $\G_\tau$.  In particular, if
$m$ is an $(\mathbb F,\mathbb P)$-martingale, the process
$(m_{t\wedge \tau}, t\geq 0)$ will be an $((\G_{t\wedge \tau},t\geq
0)$,$\mathbb Q$)-martingale (but not necessarily an $(\mathbb
G,\mathbb Q)$-martingale).
 From Theorems  \ref{Thm:G maringale characterization} and \ref{Thm:Girsanov}
 and   Corollary \ref{Cor:Girsanovimmersion}, one gets the obvious proposition.
\begin{Pro} Let
$(\Omega,\proba,\mathbb F,\mathbb G,\tau)$ be a model satisfying
the immersion property.\\ Let $(q_t(\theta), t\geq \theta) $  be a
family of positive $(\mathbb F, \mathbb P)$-martingales such that
$q_\theta (\theta)=1$ and let $\mathbb Q$ be the probability
measure with Radon-Nikod\'ym density equal to
 the $(\mathbb G,\mathbb P)$-martingale \begin{equation} Q_t^{\mathbb
G}=\indic_{\{\tau>t\}}+q_t(\tau)\indic_{\{\tau\leq t\}}\,.
\end{equation}
Then, $\mathbb Q $ and $\mathbb P$ coincide on $\G_\tau$ and the
$\mathbb P$ and $\mathbb Q$ intensities of $\tau$ are the same.
\\
Furthermore,   if $S^{\mathbb Q} $ is the $\mathbb Q$-survival
process,  the processes $( {S_t}/ {S^{\mathbb Q}_t},t\geq 0) $ and
the family $ ( \alpha_t^{\mathbb Q}(\theta) \, {S_t }/ S^{\mathbb
Q}_t ,\, t\geq \theta) $   are $(\mathbb F, \mathbb
P)$-martingales.
\end{Pro}
\proof  The first part is a direct consequence of the previous
results. It remains to note that   $Q_t^{\mathbb
F}=\frac{S_t}{S^{\mathbb Q}_t} $ and, for $t\geq \theta$,  $
\alpha_t^{\mathbb Q} (\theta)= \alpha_\theta(\theta)
\frac{q_t(\theta)}{Q_t^{\mathbb F} }= \alpha_\theta(\theta)
\frac{q_t(\theta)}{S_t }S_t^{\mathbb Q} $; hence the martingale
properties follow  from the ones  of $Q_t^{\mathbb F}$ and
$q_t(\theta)$ .\finproof
This result admits a converse.  For the sake of simplicity we
assume   the condition $S^*_t>0, \forall t\in \R_+$. This
assumption  can be removed, using the  terminal time $
\zeta^{*,\mathbb F}$.
\begin{Pro} \label{Pro:parameters} Let
$(\Omega,\proba,\mathbb F,\mathbb G,\tau)$ be a model satisfying
the H-hypothesis, with the decreasing survival process
$S_t=\exp(-\int_0^t\lambda_s^{\mathbb F}\eta(ds))$.\\
Let   $(   \alpha_t^*(\theta),t\geq\theta)$ be a given family,
where, for all $\theta >0$, $ \alpha ^*(\theta)$ is a non-negative
process and define $  S_t^*=1-\int_0^t
\alpha^*_t(\theta)\eta(d\theta)$. Assume   that $ S^*_\infty=0$
and $  S^*_t
>0, \forall t\in \R_+$ and that \begin{equation}\label{conditions2}
\left\{
\begin{array}{l}\forall \theta,\,  \alpha^* _ \theta(\theta) =
 S^*_\theta \lambda^{\mathbb F}_\theta=\alpha_ \theta(\theta) \frac{S^*_\theta}{S_\theta}\;
 \\[2mm]
 \text{the~processes~}\big(\frac{S_t}{  S^*_t}, t\geq 0
\big)\;\;\text{and} \;\;\big( \alpha^*_t(\theta)\frac{S_t}{
S^*_t},t\geq \theta \big)\; \text{are}\;\; (\mathbb F,
\proba)\text{-martingales}.
\end{array}\right.
\end{equation}
\\Let \begin{equation}\label{rniko}Q_t^{\mathbb G}:=\indic_{\{\tau>t\}}+\frac{
\alpha^*_t(\tau)  } {\alpha_{\tau}(\tau)  }\, \frac{ S_t } { S^*_t
}\indic_{\{\tau\leq t\}},
\end{equation} and $\mathbb Q$ be the probability measure with
Radon-Nikod\`ym density  the $(\mathbb G, \proba)$-martingale
$Q^{\mathbb G}$. Then,   $\Q$ is equal to $\proba$ on $\G_{\tau}$
and
\begin{equation}\label{conditions}\lambda ^{\Q,\mathbb
F}= \lambda ^{\mathbb F}, \quad\alpha_t^{\Q}(\theta)=\
\alpha^*_t(\theta)\,,\, \forall
  t \geq \theta \quad\text{and}\quad S ^{\mathbb Q}= S^*
\end{equation}
\end{Pro}
\proof    We set $$ q _t(\theta) = \frac{  \alpha^* _
t(\theta)}{\alpha _\theta(\theta)}    \frac{S_t }{
 S^*_t}. $$
 Note that $q_t(t)=1$ since $   \alpha^* _
s(s)         =
 S^*_s \lambda ^{\mathbb F}_s=   S^*_s   \alpha _s(s)/S_s$.
For every $\theta $, the processes $(q_t(\theta), t\geq \theta) $
are martingales since
 $(\alpha^*_ t(\theta) S_t /S^*_t,\, t\geq \theta)$ are martingales.
  From Theorem \ref{Thm:Girsanov}, the $\mathbb F$-projection of the Radon-Nikod\'ym density
  $Q^{\mathbb G}$ is
  $$Q^{\mathbb F}_t=
 S_t +\int_0^t \frac{  \alpha^*
 _t(s)}{\alpha_t(s)}\frac{S_t}{  S^*_t}\alpha _t(s)  \eta(ds) =S_t \left(1 +\int_0^t    \alpha^*
 _t(s) \frac{1}{  S^*_t} \eta(ds)\right),$$
 and  the survival probability is $$\frac{S_t}{Q^{\mathbb F}_t}= \big(1   +\frac{1}{
S^*_t} \int_0^t    \alpha^*
 _t(s) \eta(ds)\big) ^{-1}=  S^*_t \big(  S^*_t   +  \int_0^t    \alpha^*
 _t(s) \eta(ds)\big) ^{-1}=  S^*_t. $$
It remains to note that the condition $   \alpha^* _ s(s) =
 S^*_s \lambda ^{\mathbb F}_s$ is equivalent to the fact that the
 intensity  of $\tau$ under $\mathbb Q$ is $\lambda ^{\mathbb F}$.
 \finproof

\section{Conclusion} Our study relies on the impact
of information related to the default time on the market.

Starting from a default-free model, where some assets are traded
with the knowledge of a reference  filtration $\ff$, we consider
the case where the participants of the market take into account
the possibility of a default in view of trading default-sensitive
asset. If we are only concerned by what happens up to the default
time, the natural assumption is to assume immersion property with
stochastic intensity process adapted to the default-free market
evolution.\\
The final step is to anticipate that the default  should have  a
large impact on the market, as now after the crisis. In
particular, with the non constant ``after default" density, we
express how the default-free market is modified after the default.
In addition, hedging strategies of default-free
contingent claims are not the same in the both universes.\\

In a following paper \cite{EJJ2009}, we shall apply this
methodology to several default times, making this tool  powerful
for correlation of defaults. In another paper, we shall provide
explicit examples of density processes, and give some general
construction of these processes.


\end{document}